%This is a plain tex paper. 
\magnification=\magstep1
\vbadness=10000
\hbadness=10000
\tolerance=10000

\def\C{{\bf C}}   % complex numbers
   % rational numbers
\def\R{{\bf R}}   % real numbers
\def\Z{{\bf Z}}   % integers

\beginsection What is moonshine?

Richard E. Borcherds, 
\footnote{$^*$}{ Supported by a Royal Society
professorship.}

D.P.M.M.S.,
16 Mill Lane, 
Cambridge, 
CB2 1SB,
England.

e-mail: reb@dpmms.cam.ac.uk

home page:  www.dpmms.cam.ac.uk/\hbox{\~{}}reb

\bigskip
This is an informal write up of my talk at the I.C.M. in Berlin. It gives some
background to Goddard's talk [Go] about the moonshine conjectures. 
For other survey talks about similar topics see [B94], [B98],
[LZ], [J], [Ge], [Y].

The classification of finite simple groups shows that every finite simple group
either fits into one of about 20 infinite families, or is one
of 26 exceptions, called sporadic simple groups. 
The monster simple group is the largest of the sporadic finite simple groups,
and was discovered by Fischer and Griess [G].
Its order is 
$$
\eqalign{
&8080,17424,79451,28758,86459,90496,17107,57005,75436,80000,00000 \cr
=& 
2^{46}.3^{20}.5^9.7^6.11^2.13^3.17.19.23.29.31.41.47.59.71\cr
}$$
(which is roughly the number of elementary particles in the earth). 
The smallest irreducible representations have dimensions
$1, 196883, 21296876, \ldots$.
The elliptic modular function $j(\tau)$ has the power
series expansion 
$$j(\tau) = q^{-1}+744 + 196884q + 21493760q^2+\ldots$$ where
$q=e^{2\pi i \tau}$, and is in some sense the simplest nonconstant
function satisfying the functional equations
$j(\tau)=j(\tau+1)=j(-1/\tau)$.  John McKay noticed some rather weird
relations between coefficients of the elliptic modular function and
the representations of the monster as follows:
$$
\eqalign{
1&=1\cr
196884&=196883+1\cr
21493760&= 21296876+196883+1\cr
}
$$
where the numbers on the left are coefficients of $j(\tau)$ and the
numbers on the right are dimensions of irreducible representations of
the monster. At the time he discovered these relations, several people
thought it so unlikely that there could be a relation between the
monster and the elliptic modular function that they politely told
McKay that he was talking nonsense. The term ``monstrous moonshine''
(coined by Conway) refers to various extensions of McKay's observation,
and in particular to relations between sporadic simple groups and
modular functions.

For the benefit of readers who are not native English speakers, I had
better point out that ``moonshine'' is not a poetic terms referring to
light from the moon. It means foolish or crazy ideas. (Quatsch in
German.)  A typical example of its use is the following quotation from
E. Rutherford (the discoverer of the nucleus of the atom): ``The
energy produced by the breaking down of the atom is a very poor kind
of thing. Anyone who expects a source of power from the
transformations of these atoms is talking moonshine.'' (Moonshine is
also a name for corn whiskey, especially if it has been smuggled or
distilled illegally.)

We recall the definition of the elliptic modular function $j(\tau)$. 
The group $SL_2(\Z)$ acts on the upper half plane $H$
by 
$$\pmatrix{a&b\cr c&d\cr}(\tau)={a\tau+b\over c\tau+d}.$$
A modular function (of level 1) is a function $f$ on $H$ such 
that $f((a\tau+b)/(c\tau+d))=f(\tau)$ for all ${ab\choose cd}\in SL_2(\Z)$. 
It is sufficient to assume that $f$ is invariant under the generators
$\tau\mapsto \tau+1$ and $\tau\mapsto -1/\tau$ of $SL_2(\Z)$. 
The elliptic modular function $j$ is the simplest nonconstant example, in the sense that
any other modular function can be written as a function of $j$. 
It can be defined as follows:
$$
\eqalign{
j(\tau)&= {E_4(\tau)^3\over \Delta(\tau)}\cr
&= q^{-1}+744+196884q+21493760q^2+\cdots\cr
E_4(\tau)&= 1+240\sum_{n>0}\sigma_3(n)q^n\qquad\qquad 
(\sigma_3(n)=\sum_{d|n}d^3)\cr
&= 1+240q+2160q^2+\cdots\cr
\Delta(\tau)&= q\prod_{n>0}(1-q^n)^{24}\cr
&= q-24q+252q^2+\cdots\cr }$$ A modular form of weight $k$ is a
holomorphic function $f(\tau)=\sum_{n\ge 0}c(n)q^n$ on the upper half
plane satisfying the functional equation
$f((a\tau+b)/(c\tau+d))=(c\tau+d)^kf(\tau)$ for all ${ab\choose cd}\in
SL_2(\Z)$.  The function $E_4(\tau)$ is an Eisenstein series and is a modular
form of weight 4, while $\Delta(\tau)$ is a modular form of weight 12.

The function $j(\tau)$ is an isomorphism from the quotient
$SL_2(\Z)\backslash H$ to $\C$, and is uniquely defined by this up to
multiplication by a constant or addition of a constant. In particular any
other modular function is a function of $j$, so $j$ is in some sense
the simplest nonconstant modular function.

An amusing property of $j$ (which so far seems to have no relation
with moonshine) is that $j(\tau)$ is an algebraic integer whenever
$\tau$ is an imaginary quadratic irrational number. A well known
consequence of this is that
$$\exp(\pi\sqrt{163})=262537412640768743.99999999999925\ldots$$
is very nearly an integer. 
The explanation of this is that $j((1+i\sqrt{163})/2)$ is exactly 
the integer $-262537412640768000=-2^{18}3^35^323^329^3$, and
$$
\eqalign{
j((1+i\sqrt{163})/2)&= q^{-1}+744+196884q+\cdots\cr
&=-e^{\pi\sqrt{163}}+744+(\hbox{something very small}).\cr
}$$

McKay and Thompson suggested that there should be a graded representation 
$V=\oplus_{n\in \Z} V_n$ of the monster, such that $\dim(V_n)=c(n-1)$, 
where $j(\tau)-744 =\sum_nc(n)q^n=q^{-1}+196884q+\cdots$. Obviously 
this is a vacuous statement if interpreted literally, as we could for example
just take each $V_n$ to be a trivial representation. To characterize $V$, 
Thompson suggested looking at the McKay-Thompson series
$$T_g(\tau)=\sum_nTr(g|V_n)q^{n-1}$$ for each element $g$ of the
monster. For example, $T_1(\tau)$ should be the elliptic modular
function. Conway and Norton [C-N] calculated the first few terms of
each McKay-Thomson series by making a reasonable guess for the
decomposition of the first few $V_n$'s into irreducible
representations of the monster. They discovered the astonishing fact
that all the McKay-Thomson series appeared to be Hauptmoduls for
certain genus 0 subgroups of $SL_2(\Z)$. (A Hauptmodul for a subgroup
$\Gamma$ is an isomorphism from $\Gamma\backslash H$ to $\C$,
normalized so that its Fourier series expansion starts off
$q^{-1}+O(1)$.)

As an example of some Hauptmoduls of elements of the monster, we will look at
the elements of order 2. There are 2 conjugacy classes of elements of order 2, 
usually called the elements of types $2A$ and $2B$. The corresponding 
McKay-Thompson series start off
$$\eqalign{
T_{2B}(\tau)&= q^{-1} +276 q -2048q^2+\cdots \qquad 
\hbox{Hauptmodul for } \Gamma_0(2)
\cr
T_{2A}(\tau)&= q^{-1} +4372 q +96256q^2+\cdots \qquad 
\hbox{Hauptmodul for }\Gamma_0(2)+
\cr
}$$

The group $ \Gamma_0(2)$ is $ \{ {ab\choose cd} 
\in SL_2(\Z)|c\hbox{ is even} \}$, and
the group $\Gamma_0(2)+$ is the normalizer of $\Gamma_0(2)$ in
$SL_2(\R)$. Ogg had earlier commented on the fact that the full
normalizer $\Gamma_0(p)+$ of $\Gamma_0(p)$ for $p$ prime is a genus 0
group if and only if $p$ is one of the primes 2, 3, 5, 7, 11, 13, 17,
19, 23, 29, 31, 41, 47, 59, or 71 dividing the order of the monster.

Conway and Norton's conjectures were soon proved by O. L. Atkin,
P. Fong, and S. D. Smith. The point is that to prove something is a
virtual character of a finite group it is only necessary to prove a finite
number of congruences. In the case of the moonshine module $V$,
proving the existence of an infinite dimensional representation of the
monster whose McKay-Thompson series are give Hauptmoduls requires
checking a finite number of congruences and positivity conditions
for modular functions,
which can be done by computer.

This does not give an explicit construction of $V$, or an explanation
about why the conjectures are true. Frenkel, Lepowsky,
and Meurman managed to find an explicit construction of a monster
representation $V=\oplus V_n$, such that $\dim(V_n)=c(n-1)$, and this
module had the advantage that it came with some extra algebraic
structure preserved by the monster. However it was not obvious that
$V$ satisfied the Conway-Norton conjectures. So the main problem in
moonshine was to show that the monster modules constructed by Frenkel,
Lepowsky and Meurman on the one hand, and by Atkin, Fong, and Smith on
the other hand, were in fact the same representation of the monster.

Peter Goddard [Go] has given a description of the proof of this in his
talk in this volume, so I will only give a quick sketch of this. The main steps of the proof are as follows:
\item{1.} The module $V$ constructed by Frenkel, Lepowsky, and Meurman 
has an algebraic structure making it into a ``vertex algebra''. A
detailed proof of this is given in [F-L-M].
\item{2.} Use the vertex algebra structure on $V$ and the Goddard-Thorn
no-ghost theorem [G-T] from string theory to construct a Lie algebra acted
on by the monster, called
the monster Lie algebra.
\item{3.} The monster Lie algebra is a ``generalized Kac-Moody algebra'' 
([K90]); 
use the (twisted) Weyl-Kac denominator formula to show that $T_g(\tau)$
is a ``completely replicable function''. 
\item{4.} 
Y. Martin [M], C. Cummins, and T. Gannon [C-G] proved several
theorems showing that completely replicable functions were modular
functions of Hauptmoduls for genus 0 groups.  By using these theorems
it follows that $T_g$ is a Hauptmodul for a genus 0 subgroup of
$SL_2(\Z)$, and hence $V$ satisfies the moonshine conjectures. (The
original proof used an earlier result by Koike [Ko] showing that the
appropriate Hauptmoduls were completely replicable, together with a boring
case by case check and the
fact that a completely replicable function is characterized by its
first few coefficients.)

We will now give a brief description of some of the terms above, 
starting with vertex algebras. The best reference for finding out more 
about vertex algebras is Kac's book [K]. In this paragraph we give a rather vague description. Suppose that $V$ is a commutative ring acted on by a group $G$. We can form expressions like 
$$u(x)v(y)w(z)$$ where $u,v,w\in V$ and $x,y,z\in G$, and the action
of $x\in G$ on $u\in V$ is denoted rather confusingly by $u(x)$. (This
is not a misprint for $x(u)$; the reason for this strange notation is
to make the formulas compatible with those in quantum field theory,
where $u$ would be a quantum field and $x$ a point of space-time.)
For each fixed $u,v,\ldots\in V$, we can think of $u(x)v(y)\cdots$ as
a function from $G^n$ to $V$. We can rewrite the axioms for a
commutative ring acted on by $G$ in terms of these functions.  We can
now think of a vertex algebra roughly as follows: we are given lots of
functions from $G^n$ to $V$ satisfying the axioms mentioned above,
with the difference that these functions are allowed to have
certain sorts of singularities. In other words a vertex algebra is
a sort of commutative ring acted on by a group $G$, except that 
the multiplication is not defined everywhere but has singularities. 
In particular we cannot recover an underlying ring by defining 
the product of $u$ and $v$ to be $u(0)v(0)$, because the function
$u(x)v(y)$ might happen to have a singularity at $u=v=0$. 

It is easy to write down examples of vertex algebras:
any commutative ring acted on by a group $G$ is an example. (Actually
this is not quite correct: for technical reasons we should use a
formal group $G$ instead of a group $G$.) Conversely any vertex
algebra ``without singularities'' can be constructed in this way.
Unfortunately there are no easy examples of vertex algebras that are
not really commutative rings. One reason for this is that nontrivial
vertex algebras must be infinite dimensional; the point is that if a
vertex algebra has a nontrivial singularity, then by differentiating
it we can make the singularity worse and worse, so we must have an
infinite dimensional space of singularities. This is only possible if
the vertex algebra is infinite dimensional. However there are plenty 
of important infinite dimensional examples; see for example Kac's book 
for a construction of the most important examples, and [FLM] for
a construction of the monster vertex algebra. 

Next we give a brief description of generalized Kac-Moody algebras. 
The best way to think of these is as infinite dimensional Lie algebras 
which have most of the good properties of finite dimensional reductive
Lie algebras. Consider a typical finite dimensional reductive Lie algebra $G$, 
(for example the Lie algebra $G=M_n(\R)$ of $n\times n$ real matrices). 
This has the following properties:
\item{1.} 
$G$ has an invariant symmetric bilinear form $(,)$ (for example 
$(a,b)=-Tr(a,b)$). 
\item{2.} 
$G$ has a (Cartan) involution $\omega$ 
(for example, $\omega(a)=-a^t$).
\item{3.} 
$G$ is graded as $G=\oplus_{n\in \Z} G_n$ with $G_n$ finite
dimensional and with $\omega$ acting as $-1$ on the ``Cartan
subalgebra'' $G_0$.  (For example, we could put the basis element
$e_{i,j}$ of $M_n(\R)$ in $G_{i-j}$.)
\item{4.} 
$(a,\omega(a))>0$ if $g\in G_n$, $g\ne 0$. 

Conversely any Lie algebra satisfying the conditions above is
essentially a sum of finite dimensional and affine Lie
algebras. Generalized Kac-Moody algebras are defined by the same
conditions with one small change: we replace condition 4 by
\item{4'.} 
$(a,\omega(a))>0$ if $g\in G_n$, $g\ne 0$ and $n\ne 0$. 

This has the effect of allowing an enormous number of new examples,
such as all Kac-Moody algebras and the Heisenberg Lie algebra (which
behaves like a sort of degenerate affine Lie algebra). Generalized
Kac-Moody algebras have many of the properties of finite dimensional
semisimple Lie algebras, and in particular they have an analogue of
the Weyl character formula for some of their representations, and an
analogue of the Weyl denominator formula. An example of the Weyl-Kac
denominator formula for the algebra $G=SL_2[z,z^{-1}]$ is
$$\prod_{n>0}(1-q^{2n})(1-q^{2n-1}z)(1-q^{2n-1}z^{-1}) = \sum_{n\in
\Z} (-1)^nq^{n^2}z^n.$$ This is the Jacobi triple product identity,
and is also the Macdonald identity for the affine Lie root system
corresponding to $A_1$.

Dyson described Macdonald's discovery of the Macdonald identities in [D].
Dyson found identities for $\eta(\tau)^m=q^{m/24}\prod_{n>0}(1-q^n)^{m}$ for the following values of $m$:
$$3, 8, 10, 14, 15, 21, 24, 26, 28, \ldots$$ and wondered where this
strange sequence of numbers came from.  (The case $m=3$ is just the
Jacobi triple product identity with $z=1$.) Macdonald found his
identities corresponding to affine root systems, which gave an
explanation for the sequence above: with one exception, the numbers
are the dimensions of simple finite dimensional complex Lie algebras.
The exception is the number 26 (found by Atkin), which as far as I
know has not been explained in terms of Lie algebras. It seems
possible that it is somehow related to the fake monster Lie algebra
and the special dimension 26 in string theory.

Next we give a quick explanation of ``completely replicable'' functions.
A function is called completely replicable if its coefficients satisfy 
certain relations. As an example of a completely replicable function, 
we will look at the elliptic modular function $j(\tau)-744=\sum c(n)q^n$. 
This satisfies the identity 
$$j(\sigma)-j(\tau)=p^{-1}\prod_{m>0\atop n\in
\Z}(1-p^mq^n)^{c(mn)}$$ where $p=e^{2\pi i \sigma}$, $q=e^{2\pi i
\tau}$. (This formula was proved independently in the 80's by Koike,
Norton, and Zagier, none of whom seem to have published their proofs.)
Comparing coefficients of $p^mq^n$ on both sides gives many relations
between the coefficients of $j$ whenever we have a solution of $m_1n_1=m_2n_2$ 
in positive integers, which are more or less the relations
needed to show that $j$ is completely replicable. 
For example, from the relation $2\times 2=1\times 4$ we get the relation
$$c(4)=c(3)+{c(1)^2-c(1)\over 2}$$
or equivalently
$$20245856256= 864299970+{196884^2-196884\over 2}.$$

In the rest of this paper we will discuss various extensions of the
original moonshine conjectures, some of which are still unproved. The
first are Norton's ``generalized moonshine'' conjectures [N]. If we look
at the Hauptmodul $T_{2A}(\tau)=q^{-1}+4372q+\ldots$ we notice that
one of the coefficients is almost the same as the dimension 4371 of
the smallest non-trivial irreducible representation of the baby
monster simple group, and the centralizer of an element of type 2A in
the monster is a double cover of the baby monster. Similar things
happen for other elements of the monster, suggesting that for each
element $g$ of the monster there should be some sort of graded
moonshine module $V_g=\oplus_n V_{g,n}$ acted on by a central
extension of the centralizer $Z_M(g)$.  In particular we would get
series $T_{g,h}(\tau)=\sum_n Tr(h|V_{g,n})q^n$ satisfying certain
conditions. Some progress has been made on this by Dong, Li, and Mason
[D-L-M], who proved the generalized moonshine conjectures in the case
when $g$ and $h$ generate a cyclic group by reducing to the case when
$g=1$ (the ordinary moonshine conjectures).  G. H\"ohn has made some
progress in the harder case when $g$ and $h$ do not generate a cyclic
group by constructing the required modules for the baby monster (when
$g$ is of type $2A$).  It seems likely that his methods would also
work for the Fischer group $Fi_{24}$, but it is not clear how to go
further than this.  There might be some relation to elliptic
cohomology (see [Hi]for more discussion of this), as this also
involves pairs of commuting elements in a finite group and modular
forms.

The space $V_g$ mentioned above does not always have an invariant
vertex algebra structure on it. Ryba discovered that a vertex algebra
structure sometimes magically reappears when we reduce $V_g$ modulo
the prime $p$ equal to the order of $g$. In fact $V_g/pV_g$ can often
be described as the Tate cohomology group $\hat H^0(g,V)$ for a suitable
integral form $V$ of the monster vertex algebra. This gives natural
examples of vertex algebras over finite fields which do not lift
naturally to characteristic 0. (Note that most books and papers on
vertex algebras make the assumption that we work over a field of
characteristic 0; this assumption is often unnecessary and excludes
many interesting examples such as the one above.)

We will finish by describing some more of McKay's observations about
the monster, which so far are completely unexplained.  The monster has
9 conjugacy classes of elements that can be written as the product of
two involutions of type $2A$, and their orders are 1, 2, 3, 4, 5, 6,
2, 3, 4. McKay pointed out that these are exactly the numbers
appearing on an affine $E_8$ Dynkin diagram giving the linear relation
between the simple roots. They are also the degrees of the irreducible
representations of the binary icosahedral group.  A similar thing
happens for the baby monster: this time there are 5 classes of
elements that are the product of two involutions of type $2A$ and
their orders are 2, 4, 3, 2, 1. (This is connected
with the fact that the baby monster is a ``3,4-transposition group''.)
These are the numbers on an affine
$F_4$ Dynkin diagram, and if we take the ``double cover'' of an $F_4$
Dynkin diagram we get an $E_7$ Dynkin diagram. The number on an $E_7$
Dynkin diagram are 1, 1, 2, 2, 3, 3, 4, 2 which are the dimensions of the
irreducible representations of the binary octahedral group. The double
cover of the baby monster is the centralizer of an element of order 2
in the monster. Finally a similar thing happens for $Fi_{24}.2$: this
time there are 3 classes of elements that are the product of two
involutions of type $2A$ and their orders are 2, 3, 1. (This is connected
with the fact that $F_{24}.2$ is a ``3-transposition group''.)
These are the
numbers on an affine $G_2$ Dynkin diagram, and if we take the ``triple
cover'' of an $G_2$ Dynkin diagram we get an $E_6$ Dynkin diagram. The
number on an $E_6$ Dynkin diagram are 1, 1, 1, 2, 2, 2, 3, which are the
dimensions of the irreducible representations of the binary tetrahedral
group. The triple cover of $Fi_{24}.2$ is the centralizer of an
element of order 3 in the monster.

The connection between Dynkin diagrams and 3-dimensional rotation
groups is well understood (and is called the McKay correspondence),
but there is no known explanation for the connection with the monster.

\proclaim References.

\item{[B94]} R. E. Borcherds,  Simple groups and string theory, 
     First European congress of mathematics, Paris July 1992, 
     Ed. A. Joseph and others,vol. 1 p. 411-421,  Birkhauser 1994. 
\item{[B98]} R. E. Borcherds, 
Automorphic forms and Lie algebras. Current developments in mathematics
     1996, to be published by international press, 1998. Also available from
www.dpmms.cam.ac.uk/\hbox{\~{}}reb.
\item{[C-G]}
C. J. Cummins,  T. Gannon,  Modular equations and the genus zero
property of moonshine functions. Invent. Math. 129 (1997), no. 3,
413--443.
\item{[C-N]}
{J. H. Conway, S. Norton, Monstrous moonshine,
   Bull. London. Math. Soc. 11 (1979) 308-339.}
\item{[D]}  F. J. Dyson, 
 Missed opportunities. Bull. Amer. Math. Soc. 78 (1972), 635--652.
\item{[D-L-M]} 
C. Dong, H. Li, G. Mason, Modular invariance of
trace functions in orbifold theory, preprint q-alg/9703016.
\item{[F-L-M]}
{I. B. Frenkel, J. Lepowsky, A. Meurman, 
  Vertex operator algebras and the monster, Academic press 1988. }
(Also see the announcement	
  A natural representation of the Fischer-Griess monster with the
  modular function $J$ as character, Proc. Natl. Acad. Sci. USA 81
  (1984), 3256-3260.)
\item{[Ge]} R. W. Gebert, Introduction to vertex algebras, Borcherds algebras 
and the monster Lie algebra, Internat. J. Modern Physics A8 (1993)
no 31, 5441-5503.
\item{[Go]} The Work of R.E. Borcherds, this volume, preprint math/9808136.
\item{[G-T]}
{P. Goddard and C. B. Thorn, Compatibility of the dual
   Pomeron with unitarity and the absence of ghosts in the dual resonance
   model, Phys. Lett., B 40, No. 2 (1972), 235-238.}
\item{[G]} 
R. L. Griess,  The friendly giant. 
Invent. Math. 69 (1982), no. 1, 1--102.
\item{[Hi]}
F. Hirzebruch, T.  Berger, R.  Jung,  Manifolds and modular forms. 
Aspects of Mathematics, E20. Friedr. Vieweg \& Sohn,
Braunschweig, 1992.  ISBN: 3-528-06414-5 
\item{[H]}
G. H\"ohn,  Selbstduale Vertexoperatorsuperalgebren und das
Babymonster.  [Self-dual vertex-operator superalgebras and the
Baby Monster] Dissertation, Rheinische
Friedrich-Wilhelms-Universit\"at Bonn, Bonn, 1995. Bonner
Mathematische Schriften 
286. Universit\"at Bonn, Mathematisches Institut, Bonn, 1996.
\item{[J]} E. Jurisich, Generalized Kac-Moody Lie algebras, free Lie algebras
and the structure of the Monster Lie algebra. 
J. Pure Appl. Algebra 126 (1998), no. 1-3, 233--266.
\item{[K90]}
{V. G. Kac, ``Infinite dimensional Lie algebras'',
third edition,  Cambridge
University Press, 1990. (The first and second editions (Birkhauser,
Basel, 1983, and C.U.P., 1985) do not contain the material
on generalized Kac-Moody algebras.) }
\item{[Ko]}
{M. Koike, 
On Replication Formula and Hecke Operators, Nagoya
University preprint.}
\item{[LZ]} B. Lian, G. Zuckerman, New perspectives on the
BRST-algebraic structure in string theory, hepth/921107,
Communications in Mathematical Physics 154, (1993) 613-64, and
Moonshine cohomology, q-alg/ 950101 Finite groups and Vertex Operator
Algebras, RIMS publication (1995) 87-11. 
\item{[M]} 
Y. Martin, On modular invariance of completely replicable
functions. Moonshine, the Monster, and related topics (South Hadley,
MA, 1994), 263--286, Contemp. Math., 193, Amer. Math. Soc.,
Providence, RI, 1996.
\item{[N]}
S. P. Norton, 
Appendix to G. Mason's paper in The Arcata Conference on
Representations of Finite Groups (Arcata, Calif., 1986), 181--210,
Proc. Sympos. Pure Math., 47, Part 1, Amer. Math. Soc., Providence,
RI, 1987.

\bye